\documentclass{article}
\usepackage{graphicx} % Required for inserting images

\usepackage{amsfonts,amsmath,amssymb,amsthm}
\usepackage{cleveref}
\usepackage{subcaption}
\usepackage{booktabs}

\usepackage[a4paper, margin=1in]{geometry} 

\newtheorem{proposition}{Proposition}

\def\V#1{{\mathbf #1}}
\def\P{{\mathbb P}}

\title{Pivot probabilities and norm effects in Gaussian elimination for $\beta$-ensembles}
\author{Kenji Gunawan \and John Peca-Medlin}
\date{}

\begin{document}

\maketitle

\begin{abstract}
    We analyze pivot probabilities in Gaussian elimination with partial pivoting (GEPP) for $2 \times 2$ random matrix ensembles. For GUE matrices, we resolve a previously reported discrepancy between theoretical predictions and empirical observations by deriving the exact pivot probability under standard LAPACK-style implementations. We further show that Dumitriu–Edelman tridiagonal $\beta$-ensembles agree with the earlier theoretical expectations.
\end{abstract}

\section{Introduction and background}

Gaussian elimination (GE) with partial pivoting (GEPP) is the most widely used solver for dense linear systems $A \V x = \V b$ for $A \in \mathbb C^{n \times n}$. For instance, GEPP is the default method if using the standard `$\backslash$' operator to solve this system in MATLAB, and  thus remains a staple in introductory linear algebra courses. GEPP iteratively transforms the system into upper triangular form via rank-1 updates, with row swaps applied if needed before each elimination step to ensure the pivot is maximal in magnitude in the leading column. GEPP results in the matrix factorization $PA = LU$ in $\mathcal O(n^3)$ FLOPs, where $L$ is lower triangular, $U$ upper triangular, and $P$ a permutation matrix that encodes each of the row swaps used throughout. 

In \cite{P24}, the second author studied the question of how many GEPP row swaps are needed for random matrices.\footnote{This is equivalent to studying the number of cycles in the disjoint cycle decomposition of the permutation derived from the $P$ factor (see \cite{P24}, with \cite{PZ_2024} that further highlights this connection).} For a simpler focus on a $2 \times 2$ matrix $A$, then either 0 or 1 row swaps are possible, and so this answer determines a Bernoulli random variable with success probability $p = \P(|A_{21}| > |A_{11}|)$. 

If $A$ has independent and identically distributed (iid)  standard normal entries, $N(0,1)$, then $p = \frac12$ as the resulting permutation is necessarily uniform (cf. \cite[Theorem 1]{P24}). So this question is more interesting for non-iid Gaussian models, such as the Gaussian orthogonal ensemble (GOE) or Gaussian unitary ensemble (GUE). Standard sampling methods for each involves first forming $G$ with iid $N(0,1)$ (resp. $N_{\mathbb C}(0,1)$) entries, and then forming \begin{equation}\label{eq: goe gue}
    A = \frac1{\sqrt 2}(G + G^*). 
\end{equation}
In this case, then $A_{11} \sim N(0,2)$ and $A_{21} \sim N(0,1)$ for GOE, while $A_{11} \sim N(0,1)$ and $A_{21} \sim N_{\mathbb C}(0,1)$ for GUE. A similar construction carries over to sample from the Gaussian symplectic ensemble (GSE), using instead iid standard quaternion normal entries, $N_{\mathbb Q}(0,1)$, to first form $G$ (and the Hermitian transpose is replaced by the quaternion conjugate transpose). 

The GOE, GUE, and GSE matrices are also called $\beta$-Hermite ensembles, where $\beta = 1,2,4$, respectively. These particular $\beta$-ensembles arise naturally in physics and random matrix theory applications (e.g., see \cite{Mehta}). This $\beta$ parameter explicitly manifests in terms of the eigenvalue problem for each ensemble, where the $\beta$-Hermite ensembles each have joint density 
\begin{equation}
    f_\beta(\boldsymbol \lambda) = c_\beta \prod_{i < j} | \lambda_i - \lambda_j|^\beta \exp\left(- \sum_{j = 1}^n \lambda_j^2/2\right), \qquad  c_\beta = \frac1{(2\pi)^{n/2}} \prod_{j=1}^n \frac{\Gamma\left(1 + \frac\beta 2\right)}{\Gamma\left(1 + \frac\beta 2 j\right)}.
\end{equation}
While the standard $\beta = 1,2,4$ models then have corresponding explicit matrix constructions (outlined above) that match these eigenvalue distributions, general $\beta$-ensembles did not have a corresponding matrix model until Dumitriu and Edelman provided such a construction of tridiagonal real-valued random ensembles in \cite{DE02} for each $\beta > 0$. Their construction noted first the use of Householder reflectors to tridiagonalize the Hermitian ensembles yielded matrices of the form
\begin{equation}
    H_\beta \sim  \frac1{\sqrt 2} \begin{bmatrix}
        N(0,2) & \chi_{\beta(n-1)}\\
        \chi_{\beta(n-1)} & N(0,2) & \chi_{\beta(n-2)}\\
        & \ddots & \ddots & \ddots\\
        && \chi_{2\beta} & N(0,2) & \chi_{\beta}\\
        &&& \chi_\beta & N(0,2)
    \end{bmatrix},
\end{equation}
and this tridiagonal form holds for any $\beta > 0$. We re-emphasize that these random matrix models are real-valued, unlike as was seen for GUE and GSE.

\section{Pivot movements needed on random matrices}

\begin{figure}
    \centering
    \includegraphics[width=0.7\linewidth]{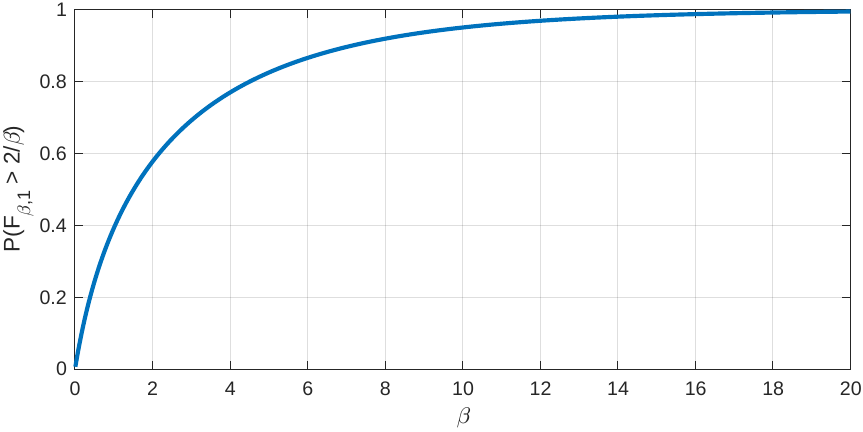}
    \caption{Plot of $2\times 2$ pivot probabilities for $H_\beta$}
    \label{fig:beta pivots}
\end{figure}

Returning to the question of the number of GEPP pivot movements needed, this can  be addressed directly for $A \sim H_\beta$ matrices of size $2 \times 2$.\footnote{Although tridiagonal $2\times 2$ matrices are then fully dense matrices, we will maintain the convention of referring to these as tridiagonal.} In this case, as $A_{11} \sim N(0,1)$ and $A_{21} \sim \chi_\beta/\sqrt 2$, then the number of GEPP row swaps needed can be determined explicitly using the $F$ distribution, where $F_{\mu,\nu} \sim (\chi^2_\mu/\mu)/(\chi^2_\nu/\nu)$ has $\mu$ numerator and $\nu$ denominator degrees of freedom, and has cumulative distribution function given by 
\begin{equation}
    \P(F_{\mu,\nu} \le x) = I_{\mu x/(\mu x + \nu)}\left(\frac\mu2,\frac\nu2\right),
\end{equation}
where $I_z(a,b)$ denotes the regularized incomplete beta function. Hence, $H_\beta$ has number of row swaps determined by a Bernoulli random variable with success probability
\begin{align}
    p_\beta^{(2)} & = \P(|A_{21}| > |A_{11}|) = \P(F_{\beta,1} > 2/\beta) = 1 - I_{2/3} \left(\frac\beta2,\frac12\right).
\end{align}
In particular, then the real-valued tridiagonal $\beta$-Hermite ensembles have pivot probability $p_1^{(2)} = 1 - \frac2\pi \arctan \sqrt 2 \approx 0.3918$, $p_2^{(2)} = \frac1{\sqrt 3} \approx 0.57735$, and $p_4^{(2)} = \frac4{3 \sqrt 3} \approx 0.76980$, respectively, for $\beta = 1,2,4$. \Cref{fig:beta pivots} shows plots of pivot probabilities for general tridiagonal $\beta$-Hermite ensembles.

\section{Discrepancy for GUE pivot comparisons}

In \cite{P24}, the second author included derivations for the pivot probabilities for GOE and GUE for the $2\times 2$ case, using a sampling distribution as given by \eqref{eq: goe gue}. For GOE, where $A_{11} \sim N(0,2)$ and $A_{21} \sim N(0,1)$ are independent, then again $\P(|A_{21}| > |A_{11}|) = \P(F_{1,1} > 2) = p_1^{(2)}$ matches exactly the computation for the $H_1$ tridiagonal construction. For GUE, however, the second author showed this computation again matched the real-valued $H_2$ computation above as now $A_{11} \sim N(0,1)$ and $A_{21} \sim N_{\mathbb C}(0,1)$; but this implicitly assumed an implementation of GEPP that used the standard complex modulus $|x + iy| = \sqrt{x^2 + y^2}$ for pivot comparisons that aligns with the $L^2$-norm when viewed as a vector in $\mathbb R^2$. Under this assumption, indeed the $H_2$ pivot probability of $p_2^{(2)} = \frac1{\sqrt 3} \approx 0.57735$ again arises. However, in a later work \cite{PZ_2024} (and as noted in a footnote there), this did not match empirical trials for this using MATLAB. 

This discrepancy between theoretical and empirical findings {stems from} standard implementations of GEPP in LAPACK \cite{lapack99}, which follow the \texttt{ZGETRF} and \texttt{DGETRF} logic—also employed by MATLAB and NumPy/SciPy in Python. {To prioritize speed and avoid the cost of computing square roots}, LAPACK compares complex-valued pivot candidates using the $L^1$ norm $\|x + iy\|_1 = |x| + |y|$ instead of the standard complex modulus. (Note that this does not impact the $H_\beta$ ensembles outlined above, as these are strictly real-valued.) The second author noted this discrepancy in \cite{PZ_2024}, and how the empirical results highly suggested the pivot probability was 2/3. We summarize this as follows:

\begin{proposition}\label{prop}
    Let $A$ be a $2 \times 2$ GUE matrix and GEPP pivot comparisons use the $L^1$ norm (as in LAPACK). Then the pivot probability is
    \begin{equation}\label{eq: pivot conundrum}
    \P(\|A_{21}\|_1 > |A_{11}|) = \P(|Z_1| + |Z_2| > \sqrt 2 |Z_3|) = \frac23
\end{equation}
for $Z_i \sim N(0,1)$ iid.
\end{proposition}
  After several conversations with  collaborators and friends across multiple institutions, no one seemed to have quite the resolve--or perhaps the patience--to tackle this problem. Instead, one such friend posed this question to a class of undergraduates in an upper division probability class; and that is how the first author enters the picture to provide the derivation of \Cref{prop}.

\begin{figure}[t]
    \centering
    \begin{subfigure}[t]{0.35\linewidth}
        \centering
        \includegraphics[width=\linewidth]{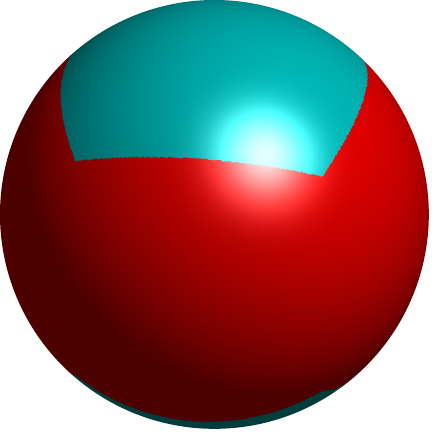}
        \caption{$|x| + |y| > \sqrt{2} |z|$}
    \end{subfigure}
    \hfill
    \begin{subfigure}[t]{0.35\linewidth}
        \centering
        \includegraphics[width=\linewidth]{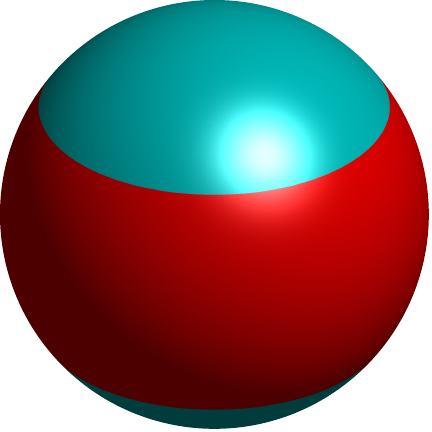}
        \caption{$\sqrt{x^2 + y^2} > \sqrt{2} |z|$}
    \end{subfigure}
    \caption{Regions on the unit sphere defined by different norms used in pivoting comparisons for $2\times 2$ GUE matrices.}
    \label{fig:pivot-regions}
\end{figure}

\section{Proof of Proposition \ref{prop}}

The approach proves to be quite straightforward, and appropriate for a graduate or advanced undergraduate probability student. A quick outline notes that \Cref{prop} can be reduced to computing the surface area for a region on the unit sphere by considering the polar decomposition of a standard iid Gaussian vector, $\V x \in \mathbb R^3$, as $\V x \sim \chi_3 \V u$ for $\V u \sim \operatorname{Unif}(S^2)$ independent of $\|\V x\|_2 \sim \chi_3$. \Cref{fig:pivot-regions} compares the corresponding regions (in red) on the unit sphere $S^2$ determined by the $|x| + |y| > \sqrt 2 |z|$, that aligns with the standard MATLAB implementation of GEPP, versus $\sqrt{x^2 + y^2} > \sqrt 2 |z|$, that would follow if the standard complex modulus was used for pivot comparisons.  Using spherical coordinates and symmetry to further reduce the computation to consider only the first octant (where $x,y,z \ge 0$), then the region $R$ of interest on the unit sphere is defined by
\begin{align}
    x + y > \sqrt 2 z &\Rightarrow \sin\varphi(\cos \theta + \sin \theta) > \sqrt 2 \cos\varphi \\
    &\Rightarrow \tan \varphi > \frac{\sqrt 2}{\cos \theta + \sin \theta} = \sec \psi
\end{align}
for $0 \le \varphi,\theta \le \frac\pi2$ and $\psi = \frac\pi4 - \theta$; we in particular note the boundary at $\varphi_0(\theta) = \arctan(\sec \psi)$, where we further note $\cos \varphi_0(\theta) = \cos\psi/\sqrt{1 + \cos^2 \psi}$. Hence, the probability involves computing the ratio of the area of this region, $A_R$, to the surface area of unit sphere in the first octant, $A_{oct} = \pi/2$. We compute then
\begin{align}
    A_R &= \int_0^{\pi/2} \int_{\varphi_0(\theta)}^{\pi/2} \sin \varphi \ d\varphi d \theta = \int_0^{\pi/2} \cos \varphi_0(\theta) \ d\theta = 2 \int_0^{\pi/4} \frac{\cos \psi}{\sqrt{1 + \cos^2 \psi}} d\psi \\&= 2\left[ \arcsin\left(\frac x{\sqrt 2}\right) \right]_0^{1/\sqrt 2} =  \frac\pi3.
\end{align}
Hence, indeed
\begin{equation}
    \P(\|A_{21}\|_1 > |A_{11}|)  = \frac{A_R}{A_{oct}} = \frac{\pi/3}{\pi/2} = \frac23.
\end{equation}

\section{Discussion}

One could instead consider the analogous question of the value of 
\begin{equation}
    p_\beta^{(1)} = \P\left( \sum_{j = 1}^\beta |Z_j| > \sqrt 2 |Z_{\beta+1}|\right)
\end{equation}
for $Z_j \sim N(0,1)$ iid. (Here, we are further focusing on integer $\beta$.) As seen above, we know the $\beta = 1,2$ probabilities evaluate, respectively, to $p_1^{(1)} = p_1^{(2)}= 1 - \frac2\pi \arctan \sqrt 2 \approx 0.3918$ and $p_2^{(1)} = \frac23 \approx 0.6667$. Each probability can similarly reduce to computing a surface area on $S^\beta$, which can be approached by similarly considering standard polar coordinates in higher dimensions (e.g., spherical coordinates in 3 dimensions); however, the resulting multi-integral is increasingly more complex. We thus leave it as an exercise for the motivated reader to derive the exact computations for these values. For reference, we include empirical results related to $10^6$ trials for each such value, that produces the following summary table for estimated values $\hat p_\beta^{(1)}$ in \Cref{tab:split-abs-normal}. 
\begin{table}[t]
    \centering
    \begin{tabular}{c c c c c c}
        \toprule
        $\beta$ & 1 & 2 & 3 & 4 & 5 \\
        \midrule
        $\hat p_\beta^{(1)}$
        & 0.3908 & 0.6658 & 0.8326 & 0.9223 & 0.9665 \\
        \midrule
        $\beta$ & 6 & 7 & 8 & 9 & 10 \\
        \midrule
        $\hat p_\beta^{(1)}$
        & 0.9864 & 0.9947 & 0.9981 & 0.9993 & 0.9998 \\
        \bottomrule
    \end{tabular}
    \caption{Estimated values of $p_\beta^{(1)}$ for $\beta = 1$ to $10$, using $10^6$ samples.}
    \label{tab:split-abs-normal}
\end{table}

Note while the $\beta = 1,2$ values correspond to actual pivot probabilities that align with standard implementations of GEPP with $2 \times 2$ GOE and GUE matrices, the $\beta = 4$ case does not align with a GSE matrix since LAPACK does not by default handle quaternion valued matrices, so a custom GEPP would need to be employed (and so could, say, assume the standard quaternion modulus $|x + \textbf{i}y + \textbf{j}u + \textbf{k}v| = \sqrt{x^2 + y^2 + u^2 + v^2}$  for pivot candidate comparisons, that would thus align with the pivot probabilities for $H_4$ of $p_4^{(2)} = \frac{4}{3 \sqrt 3} \approx 0.76980$, if so desired). 

Note by the triangle inequality the $H_\beta$ probabilities $p_\beta^{(2)}$ (cf. \Cref{fig:beta pivots}) provide lower bounds for each corresponding $p_\beta^{(1)}$ value (e.g., how $p_2^{(2)} \approx 0.57735$ is smaller than $p_2^{(1)} = \frac23$ and $p_4^{(2)} \approx 0.76980$ is smaller than the estimate $\hat p_4^{(1)} = 0.9223$, while $p_1^{(1)} = p_1^{(2)}$). For general integer $\beta$, what is the exact value of $p_\beta^{(1)}$? Closed-form expressions for higher $\beta$ appear intractable via direct integration methods  and warrant further investigation.

Pivot movements can serve as a bottleneck for GEPP, as they hinder data movements on parallel architectures \cite{Pa95}. Hence, other authors are invested in pivoting strategies that minimize communication (for example, see \cite{calu}) or preprocessing the system to avoid pivoting altogether \cite{baboulin,Pa95,PT23}. The above discussion suggests additionally the choice of norm for complex-valued pivot candidate comparisons can also be a factor in reducing communication costs.

Additionally, future work can explore general GEPP induced permutations for GOE, GUE, and GSE, as well as $H_\beta$ for $n > 2$. While these do not induce uniform permutations (see \cite{P24}; this can be justified as the first pivot selection is not uniform), while at least for GOE and GUE they do appear to be close to uniform (in a permuton sense) for sufficiently large $n$. For $H_\beta$, induced permutations are much farther from uniform due to the tridiagonal form. The limiting permutations appear to concentrate much closer to the identity permutation as $n$ grows. Future work can further explore this direction of research.

\section{Acknowledgments}

The second author thanks Michael Conroy, Adrien Peltzer, and Chenyang Zhong for early discussions on this problem, and Morris Ang and Rob Webber for helpful conversations and for sharing the question with students, which ultimately led to this collaboration.

\bibliographystyle{plain} 
\bibliography{references}

\end{document}